\documentclass[12pt,a4paper,reqno]{article}
\usepackage{amsmath,amsthm}
\usepackage{amssymb,latexsym,amsmath}
\usepackage{amscd}
\usepackage{amsfonts}
\input amssymb.sty
\chardef\bslash=`\\ 
\newtheorem{theorem}{Theorem}
\newtheorem{lemma}{Lemma}

\newtheorem{corollary}{Corollary}

\begin{document}
\title{Some remarks on autoequivalences of categories}
\author{Grigori~I. Zhitomirski\thanks{This research is partially supported by THE ISRAEL SCIENCE FOUNDATION
 founded by The Israel Academy of Sciences and Humanities - Center of
 Excellence Program.}\\
   Department of Mathematics and Statistics\\
   Bar Ilan University\\
   Ramat Gan 52900, Israel}
\date{}
\maketitle

\begin{abstract}
Prof. Boris I. Plotkin \cite{SevLect,AlgGeom} drew attention to the question when an equivalence between two
categories is isomorphic as a functor to an isomorphism between them. It turns out that it is quite important
for universal algebraical geometry and concerns mainly the categories $\Theta \sp 0 (X) $ of free universal
algebras of some variety $\Theta $ free generated by finite subsets of $X$. In the paper, a complete answer to
the Plotkin's question is given: there are no proper autoequivalences of the category $\Theta \sp 0 (X) $. Also
some connected problems are discussed.
\end{abstract}
\section*{Introduction}

Prof. Boris I. Plotkin set a question arose by studying of the universal algebraical geometry
\cite{SevLect,AlgGeom}. Without going into details, this question can be formulated in the following way. Let
$\mathcal V$ be a variety of universal algebras. Consider the category $\Theta$ whose objects are all algebras
from $\mathcal V$ and whose morphisms are all homomorphisms of them. Fix an infinite set $X$. Let $\Theta \sp
0(X) $ be the full subcategory of $\Theta $ containing only free $\mathcal V-$algebras over finite subsets of
the set $X$. The question is: if there are autoequivalences of $\Theta \sp 0 (X) $ that are not isomorphic to
any automorphism of this category? Below, we give a negative answer to this question, i. e. we show that every
autoequivalence of the category $\Theta \sp 0 (X) $ is isomorphic to an automorphism of this category, and
present some results concerning that theme. Properly to say Prof. B.I. Plotkin regards the mentioned question to
be important for the universal algebraical geometry. The reasons for this opinion are the following ones.

Let $\mathcal C$ be an arbitrary small category. Let $End{\mathcal C}$ be the monoid of all endofunctors $F:
{\mathcal C} \to {\mathcal C}$ and  $Aut {\mathcal C}$ be the group of all endofunctors that are automorphisms
of this category. The relation "two functors are isomorphic (natural equivalent)"  is a congruence relation on
the monoid $ End {\mathcal C}$. The corresponding quotient monoid is denoted by $End\sp 0 {\mathcal C}$ and the
group of all invertible elements of it is denoted by $Aut \sp {0}\ \mathcal C$. We have the monoid homomorphism:
 $\eta : End {\mathcal C} \to  End\sp 0 {\mathcal C}$ and the induced homomorphism
 $\eta \sp * : Aut {\mathcal C }\to  Aut \sp {0}\ {\mathcal C}$. The following diagram is commutative.
$$
\CD
Aut {\mathcal C} @> \eta \sp * >> Aut \sp 0 {\mathcal C} \\
@VVV                                             @VVV \\
 End {\mathcal C} @ >> \eta  >                   End\sp 0 {\mathcal C}
\endCD
$$
Here the arrows on each side are the inclusion maps.

The case the homomorphism $\eta\sp *$ being surjective is very important. It is easy to prove that this holds if
and only if every autoequivalence of the category $\mathcal C$ is isomorphic to an automorphism of this
category. Therefore we are interested to know if the category $\Theta \sp 0 (X)$ satisfies this condition.

There is one more situation that leads to the same question but this time for two categories. Let $H$ be an
object of the category $\Theta $, i.e. it is an algebra from $\mathcal V$. Let $W$ be an object of the category
$\Theta \sp 0 =\Theta \sp 0 (X)$, i.e. it is a free algebra in $\mathcal V$  over a finite subset of $X$. The
set $Hom(W,H)$ is regarded as a an affine space and its elements are called the points.  There is a Galois
correspondence between binary relations in $W$ regarded as sets of equations and the sets of points. The closed
subsets under this correspondence are called algebraical sets and $H$-closed congruences respectively. Also a
notion of algebraical variety can be introduced. Thus with every algebra $H$, some geometry is connected, and
the question appears what it means that two algebras have the same geometry. This question leads to the problem
of recognizing if the corresponding two categories merely equivalent or what is more isomorphic.

In the Section~2 we give a necessary and sufficient condition for an equivalence of two small categories to be
isomorphic to an automorphism of them (Theorem~\ref{bijection}) and obtain from it the mentioned answer to the
Plotkin's question. In the section~3 we obtain some similar results for arbitrary categories represented in the
category of sets.
\section{Preliminaries}

We use the usual notions and notation of category theory (see for example [3]). Let  $\mathcal {C}$ and
$\mathcal {D}$ be two arbitrary categories. A functor $F: \mathcal {C} \to \mathcal {D}$ is called an {\it
equivalence } if

1) it is faithful and full, i.e. for every two objects $A,B $ from $\mathcal {C}$ the functor $F$ determines
bijection of $Hom (A,B)$ onto $Hom(F(A),F(B))$;

2) every object from $\mathcal {D}$ is isomorphic to $F(A)$ for some object $A$ of $\mathcal {C}$.

These two conditions mean, that $F$ induces an isomorphism between skeletons of these two categories. An
equivalence $F: \mathcal {C} \to \mathcal {C}$ is called an {\it autoequivalence } of the category $\mathcal
{C}$.

We recall that a morphism (or a natural transformation) from a functor $F: \mathcal {C} \to \mathcal {D}$ to a
functor $G: \mathcal {C} \to \mathcal {D}$ is a map $\alpha$ which assigns to every object $A$ of $\mathcal {C}$
a morphism $\alpha {\sb A }: F(A) \to G(A)$ such that for every morphism $f: A\to B$ of $\mathcal {C}$ the
following diagram commutes:
$$
\CD
F(A)@> \alpha \sb A >> G(A) \\
@VF(f)VV              @VVG(f)V \\
F(B) @ >> \alpha \sb B >  G(B)
\endCD
$$
A functor from a category to itself is called an endofunctor. The category of all functors from a category
$\mathcal {C}$ to a category $\mathcal {D}$ is denoted by $Funct({\mathcal {C}},{\mathcal {D}} )$. The category
of endofunctors of a category $\mathcal {C}$ is denoted by $Funct\sb {\mathcal {C}}$ . Equivalences and
isomorphisms from $\mathcal {C}$ to $\mathcal {D}$ (if they exist) are objects of the category $Funct({\mathcal
{C}},{\mathcal {D}} )$, and an equivalence which is not isomorphic ( in sense of this category ) to any
isomorphism is called a proper equivalence. Every autoequivalence of $\mathcal {C}$ that is isomorphic to the
identity functor $\bold Id \sb {\mathcal {C}}$ is called an inner autoequivalence as well an automorphism (that
certainly is an autoequivalence) is called an inner automorphism if it is isomorphic to {\bf Id}.

For every object $A$ of $\mathcal {C}$, we denote by $[A]$ the class of all $\mathcal {C}$-objects that are
isomorphic to $A$. Let $Sk = Sk\sb {\mathcal C}$ be a skeleton of $\mathcal {C}$. For every object $A$ there
exists an unique common object $\overline {A}$ from $Sk$ and $[A]$.  Chose an isomorphism $u\sb A : A\to
\overline {A}$ for every object $A$ .  The next result is well known.

\begin{lemma}\label{onSc} Let  $\nu \sb {sk}$ be a map that assigns to every object $A$ the object $\overline {A}$
and to every morphism $f:A\to B$ the morphism $\overline {f} : \overline {A} \to \overline {B}$ by the following
way: $\overline {f} =u\sb A \sp {-1} f u\sb B $. Then the map $\nu \sb {sk}$ is an inner autoequivalence of the
category $\mathcal {C}$.
\end{lemma}

It is easy to see that there are categories admitting a proper autoequivalence. As an example, every category
$\mathcal {C}$ is suitable if it satisfies the condition: its skeleton $Sk = Sk\sb {\mathcal C}$ admits some
automorphism $\psi$ that does not fix at least one object but there are no non-trivial automorphisms of the
category $\mathcal {C}$ itself. Indeed, let $\nu \sb {sk}$ be the autoequivalence defined in the Lemma above.
Thus $\psi \circ \nu \sb {sk}$ is an autoequivalence but if $\psi (X)\not = X$ for some $X\in Sk$ then $\psi
\circ \nu \sb{sk}(X)$ is not isomorphic to $X$. Hence $\psi \circ \nu \sb {sk}$ is a proper autoequivalence.

Consider, for instance, four items $a,b,c,d$ as objects and pre-order relation defined by following arrows:
$a\to b, a\to c\leftrightarrows d$ and of course loops for every object. Let $\mathcal {C}$ be the corresponding
category. There is non-trivial automorphism of $Sk\sb {\mathcal C}$ but there are no non-trivial automorphisms
of $ {\mathcal C}$.

From this point of view the following fact is of interest and it answers the mentioned Plotkin's question. The
proofs of results presented in the next section are quite simple but the author has not met them in the
literature and therefore explains them explicitly.

\section{There are no proper autoequivalences of the category $\Theta \sp 0 (X) $}

\begin{theorem}\label{bijection} Let $\mathcal {C}$ and $\mathcal {D}$ be two arbitrary small categories.
An equivalence $\alpha :\mathcal C \to \mathcal D$ is isomorphic to an isomorphism if and only
if for every to objects $X$ and $Y$ such that $\alpha (X) = Y$ the sets $[X] $ and $[Y]$ are of equal
cardinality.
\end{theorem}

\begin{proof} Let $\pi :\mathcal C \to \mathcal D$ be an equivalence. Let
$\nu =\nu \sb {sk} :\mathcal {D} \to Sk \sb {\mathcal D }$. Consider $\alpha =\nu \circ \pi$. If $\pi (X)=Y\sb
1$ and $\alpha (X)=Y\sb 2$ then objects $Y\sb 1$ and $Y\sb 2$ are isomorphic and hence $[Y\sb 1 ]= [Y\sb 2]$. On
the other hand, since $\nu $ is isomorphic to the identity functor $\bold Id \sb {\mathcal D}$, an equivalence
$\pi$ is isomorphic to an isomorphism if and only if so $\nu \circ \pi$ is. These two facts mean that it is
sufficient to proof our statement for the equivalence $\alpha$.

If $\alpha$ is isomorphic to an isomorphism $\varphi$, then for every object $Y$ of $Sk \sb {\mathcal D }$,
$\varphi$ determines an one-to-one map of $\alpha \sp {-1} (Y) $ onto $[Y]$. It means that for every to objects
$X$ and $Y$ such that $\alpha (X) = Y$ the sets $[X] $ and $[Y]$ are of equal cardinality.

If this last condition is satisfied then we have an one-to-one map $\varphi \sb Y :\alpha \sp {-1} (Y) \to [Y]$
for every object $Y$ of $Sk \sb {\mathcal D }$ . So we can define for every object $A$ the object $\varphi (A)
=\varphi \sb Y (A)$ if $A \in \alpha \sp {-1} (Y)$. Let $A \in \alpha \sp {-1} (Y\sb 1), B \in \alpha \sp {-1}
(Y\sb 2)$ and $f:A\to B$ be a morphism. Then one has the morphism $\alpha (f): \alpha (A) \to \alpha (B)$. The
morphism $\varphi (f): \varphi (A) \to \varphi (B)$ can be uniquely defined from the following diagram in order
to make it commutative:
$$
\CD
\varphi (A)@> u\sb {\varphi (A)} >> Y\sb 1 \\
@V\varphi (f)VV              @VV\alpha (f)V \\
\varphi (B) @ >> u\sb {\varphi (B)} >  Y\sb 2
\endCD
$$
It is oblivious that the defined map $\varphi$ is an isomorphism of the categories $\mathcal C \to \mathcal D$
and it is isomorphic to $\alpha$.
\end{proof}

The next result follows easily from the theorem above.

\begin{theorem}\label{mainCond}

1. In order to a small category $\mathcal C$ has no proper autoequivalences it is necessary and sufficient that
for every automorphism $\gamma$ of $Sk$ and for every two objects $X$ and $Y$ of $Sk$ such that $\gamma (X)= Y$,
the sets $[X]$ and $[Y]$ are of equal cardinality. Particularly this condition is satisfied if all automorphisms
of the category $Sk\sb {\mathcal C} $ leave fixed its objects.

2. If a small category $\mathcal {C}$ satisfies the condition that for every two its objects $A,B$ the set $[A]$
and $[B]$ have the same cardinality then every autoequivalence of $\mathcal {C}$ is isomorphic to an
automorphism of this category.
\end{theorem}

Now consider the main question mentioned in Introduction. Let $\mathcal V $ be a variety of universal algebras
and $X$ be an infinite set.  The free $\mathcal V $-algebra $W=W(X\sb 0)$  corresponds to every finite subset
$X\sb 0$ of $X$ being free generated by this subset. Let $\Theta\sp 0 (X)$ be the category which objects are all
such free algebras an which morphisms are homomorphisms of $\mathcal V $-algebras. For every finite subset $X\sb
0 \subset X$ the set $[W(X\sb 0)]$ of $\Theta\sp 0 (X)$-objects has the same cardinality that the set $X$. Hence
we have

\begin{theorem}\label{free} Let $\mathcal V \sb 1$ and $\mathcal V \sb 2$ be two varieties of universal algebras. Let $X$
be an infinite set and $\Theta \sb 1 \sp 0 (X),\; \Theta \sb 2 \sp 0 (X) $ be two corresponding categories of
free algebras.  Then every equivalence of this categories is isomorphic to an isomorphism. Thus if such
categories are equivalent they are isomorphic. Further, for a categories of the kind $\Theta \sp 0 (X)$, every
its autoequivalence is isomorphic to an automorphism of this category.
\end{theorem}

Now we give more example of a category which has proper autoequivalences. Let $E$ be an ordered set such that
there is a non-trivial automorphism $\varphi $ of it. Let $e\in E$ be an element such that $\varphi (e) \neq e$.
Add to $E$ a new element $a$ and consider on the set $E'=E\cup \{a\}$ the pre-order relation $\rho$ that will be
obtained by adding two pairs $(a,e)$ and $(e,a)$ to the existing order relation on $E$. Let $\mathcal C$ be the
category defined on the base of this pre-order relation. The map $\varphi '$ that acts like $\varphi $ on $E$
and takes $a$ to $\varphi (e)$ (i.e. we have $\varphi '(a)=\varphi '(e)$) determines an proper autoequivalence
according to Theorem~~\ref{bijection} because $[e]=\{e,a\}$  and $[\varphi ' (e)] =\{\varphi ' (e)\}$.

The proof of Theorem~\ref{bijection}  uses the Axiom of Choice but it is possible to remove references to this
axiom in some simple cases, for example for the category $\Theta \sp 0 (X)$ over a countable set $X$ . One can
try to extend Theorem~\ref{bijection} up to non-small categories. But this way leads into jungle of non-sets
cardinalities. It concerns axioms which are independent from usual set-theoretic axioms. This reason explains
why it is impossible to give an example of a big category (like the category of universal algebras for some
variety) with proper autoequivalences.

Because the equality of cardinalities for every two sets of kind $[A]$ is not so easy to see we give in the next
section some other conditions for an autoequivalence to be isomorphic to an automorphism.

\section{Categories with a forgetful functor}

We consider such categories $\mathcal C$ that are represented in the category $Set$ (the category of all sets
and maps), that is, there exists a faithful functor  $Q: \mathcal C \to Set$. Such a functor is called a
forgetful functor. If $\mathcal C$ is a category of universal algebras, then the forgetful functor is usually
the natural forgetful functor, which assigns to every algebra $A$  the underlying set $\vert A \vert$ and to
every homomorphism itself as a mapping, but not only this case. Additionally we assume that the category
$\mathcal C$ satisfies the following condition for every two objects  $A$ and $B$:

(*) If   $Q(A)$ and $Q(B)$ have the same cardinality then, for every bijection $u :Q(A)\to Q(B)$, there exists
an unique object $C$ of $\mathcal {C}$ and an isomorphism $u \sp * :C\to B$ such that  $Q(C)= Q(A)$ and $Q(u \sp
*) = u$.

This requirement is a weakening of the condition to be abstract for a subcategory of a category of sets supplied
with some structure (for example, an algebraic structure). Indeed, every bijection $u :Q(A)\to Q(B)$ induces the
unique structure on the set $Q(A)$  that the map $u$ is an isomorphism of the obtained object $C$ onto $B$. Our
requirement is that this object belongs to the subcategory under consideration.

\begin{theorem}\label{generConD} Let $\pi :\mathcal {C} \to \mathcal {C}$ be an autoequivalence. If for every object
$A$ of $\mathcal {C}$ there exists a bijection $u \sb A : Q (A) \to Q\pi (A)$ such that for every morphism $f:A
\to B$ of $\mathcal {C}$ the following diagram commutes:
$$
\CD
Q(A)@> u \sb A  >> Q\pi (A) \\
@VQ(f)VV              @VVQ\pi (f)V \\
Q(B) @ >> u \sb B >  Q\pi (B)
\endCD
$$
then $\pi $ is isomorphic to an automorphism $F $ of $\mathcal {C}$.
\end{theorem}
\begin{proof}
According to hypothesis, we have for every object $A$ the bijection $u \sb A : Q (A) \to Q\pi (A)$. Because of
the property (*), there exists an unique object $F(A)$ and an unique isomorphism $u \sp * \sb A :F(A)\to \pi
(A)$ such that :
     $$Q(F(A)) =Q(A)\qquad  {\rm and}\qquad  Q(u \sp * \sb A )= u \sb A .$$
 Next, define for every morphism  $f:A \to B$ of $\mathcal {C}$ the morphism
 $F (f): F(A) \to F(B)$  by the following way: $F(f)= u \sp * \sb A \pi (f)(u \sb B \sp * )\sp {-1}$.
Thus we have the following commutative diagram:
$$
\CD
F(A)@> u \sp * \sb A  >> \pi (A) \\
@VF(f)VV              @VV\pi (f)V \\
F(B) @ >> u \sp * \sb B >  \pi (B)
\endCD
$$
It means that $F$ is an endofunctor of the category $\mathcal C $ and it is isomorphic to the given
autoequivalence $\pi$. It remains to proof that $F$ is an automorphism of $\mathcal C $.

Since $F$ is autoequivalence, we have only to show that $F$ is an bijection with respect to objects. Let $F(A)
=F(B)$. Then we have an isomorphism $g=(u \sp * \sb A )\sp {-1}u \sp * \sb B  : \pi (A) \to \pi (B)$. Thus there
exists an unique isomorphism $f: A\to B$ such that $\pi (f) =g$. According to this definition of the morphism
$f$ and hypothesis of the theorem, the following diagram commutes:
$$
\CD
Q(A)@> u \sb A  >> Q\pi (A) \\
@VQ(f)VV              @VVQ(g)V \\
Q(B) @ >> u \sb B >  Q\pi (B)
\endCD
.$$

Having in the mind that $Q(A)=Q(F(A))=Q(F(B))=Q(B)$, we obtain:
$$ Q(f)= u \sb A u \sb A \sp {-1}u  \sb B u  \sb B \sp {-1} =1\sb {Q(A)}.$$
Because of the condition (*) it means that $A=B$.

Now let $S$ be an arbitrary object of $\mathcal C$. Since $\pi$ is an autoequivalence, there exists an object
$T$ such that there exists an isomorphism $i:S\to \pi (T)$. Thus we obtain a bijection $g =Q(i)u \sb T \sp
{-1}:Q( S) \to Q(T) $. Let $V$ be an object from $\mathcal C$  such that $Q(V)=Q(S)$, $g\sp * : V\to T$ is an
isomorphism, and $Q(g\sp *)= g$. Then the following diagram commutes:
$$
\CD
Q(V)@> u \sb V  >> Q\pi (V) \\
@Vg VV              @VVQ\pi (g \sp *)V \\
Q(T) @ >> u \sb T >  Q(\pi (T))
\endCD
.$$

Let $h=u \sp * \sb V \ \pi (g \sp *)\ i \sp {-1} :F(V)\to S$. It is clear that $h$ is an isomorphism and
$Q(F(V))=Q(S)$. It follows from the diagram above and the definition of the bijection $g$ that
$$Q(h)=u \sb V \ u \sb V \sp {-1} \ g \ u \sb T \ Q(i\sp {-1}) =
Q(i)\ u \sb T \sp {-1}\ u \sb T \ Q(i\sp {-1}) =1\sb {Q(S)}$$ and hence $F(V) =S$. This completes the proof.

\end{proof}

The hypotheses of Theorem~\ref{generConD} are satisfied in many important cases.

\begin{theorem}\label{saferepr} If the functor $Q:  \mathcal C \to Set $ is  represented by an
object $W$  and if an autoequivalence $\pi$ satisfies the condition that $W$ is isomorphic to $\pi (W)$  then
$\pi$ satisfies the hypotheses of Theorem~\ref{generConD} and hence it is isomorphic to an automorphism of
$\mathcal C$.
\end{theorem}

\begin{proof} Let $W$ and $ \pi (W)$  are isomorphic according to hypotheses. Thus the functor $Q$ is
represented by $\pi (W)$ too. We have two isomorphisms in the category $Funct(\mathcal C ,Set )$, namely $\bold
g :Q \to Hom(W,-)$ and $\bold h :Q \to Hom(\pi (W), -)$. Since $\pi $ induces a bijection between sets
$Hom(A,B)$ and $Hom(\pi (A),\pi (B))$ for every pair $(A,B)$ of objects, there is a bijection $u \sb A : Q (A)
\to Q\pi (A)$ for every object $A$ defined as the composition of bijections:
$$
\CD
 Q(A)@>\bold g \sb a >>Hom (W,A)@>\pi >>Hom(\pi (W),\pi (A)) @> \bold h \sp {-1} \sb {\pi (A)}>>Q\pi (A).
\endCD
$$

The fact that the family $u \sb A : Q (A) \to Q\pi (A)$ of bijections satisfies the hypotheses of
Theorem~\ref{generConD} one can see from the following diagram in which every of three squares is commutative
for every morphism $f:A\to B$:
$$
\CD
Q(A)@>\bold g \sb a >>Hom (W,A)@>\pi >>Hom(\pi (W),\pi (A)) @> \bold h \sp {-1} \sb {\pi (A)}>>Q\pi (A) \\
@Vf VV              @VVHom(W,f)V           @VVHom(\pi (W),\pi (f))V                       @VVQ\pi (f)V \\
Q(B)@>\bold g \sb a >>Hom (W,B)@>\pi >>Hom(\pi (W),\pi (B)) @> \bold h \sp {-1} \sb {\pi (B)}>>Q\pi (B)
\endCD
.$$

\end{proof}

The first consequence from Theorem~\ref{saferepr} is the following result:

\begin{corollary}\label{set} Every full subcategory of $Set$ containing an one-element set as an object has no proper
equivalences, moreover all its autoequivalences are inner.
\end{corollary}

Let $\mathcal V$ be a variety of universal algebras and $\Theta $ be the corresponding category. Let $Q: \Theta
\to Set$ be the forgetful functor. It is clear that condition (*) of the previous section is satisfied.
Moreover, the functor $Q$ is represented by the free algebra $W$ in $\mathcal V $ over an one-element set. Thus
according to Theorem~~\ref{saferepr}, if  $W$ is isomorphic to $\pi (W)$ for an autoequivalence $\pi$ of $\Theta
$, then $\pi $ is isomorphic to an automorphism of this category. To be sure, that the same fact is valid for
every full subcategory of $\Theta $ if it satisfies the condition (*) and contains a free algebra $W$ over a
singleton.

So the question is if for an autoequivalence $\pi$ the algebras $W$ and $\pi (W)$ are isomorphic. Since $\pi$
induces an isomorphism between the monoids $END(W)$ and $END(\pi (W))$ of endomorphisms of these two algebras it
is sufficient to check if the fact that these monoids are isomorphic implies that the mentioned algebras are
isomorphic too. If it is true we can conclude that every autoequivalence is isomorphic to an automorphism of
this category. Hence we have

\begin{theorem}\label{freeAlg} Let $\mathcal C $ be the category of all free algebras in a variety $\mathcal V $ of
universal algebras. Let $W$ be the free algebra $W$ in $\mathcal V $ over an one-element set. If  any free
algebra $F$ in $\mathcal V $ is isomorphic to $W$ provided that monoids $END(W)$ and $END(F)$ are isomorphic
then the category $\mathcal C $ has no proper autoequivalence.

\end{theorem}
For many known varieties, the monogenic free algebras satisfy the condition mentioned above, in particular, the
categories of all free semigroups, of all free inverse semigroups or of all free groups have no proper
autoequivalence.

I am very thankful to Prof. Boris I. Plotkin who has called my attention to the mentioned above problems and has
discussed the results with me.


\end{document}